\documentclass[12pt]{amsart}
\usepackage{amsmath}%
\usepackage{amsfonts}%
\usepackage{amssymb}%
\usepackage{graphicx}
\usepackage{tikz-cd}
\usepackage{enumerate}
\usepackage{rotating}
\usepackage{rotfloat}
\usepackage{caption}
\usepackage[all]{xy}

\newtheorem{teo}{Theorem}[section]
\newtheorem{lemma}[teo]{Lemma}
\newtheorem{cor}[teo]{Corollary}
\newtheorem{prop}[teo]{Proposition}
\newtheorem{defi}[teo]{Definition}
\newtheorem{example}[teo]{Example}

\newtheorem{conj}[teo]{Conjecture}
\newtheorem{prob}[teo]{Problem}

\floatstyle{plain}
\floatname{diagram}{Diagram}
\newfloat{diagram}{tbp}{lop}[section]

\begin{document}

\title[Some results about Colored Tverberg Theorem]{Some results about Colored Tverberg Theorem}


\author[\hspace{1cm}Carlos H.F. Poncio]{Carlos H. F. Poncio}
\address{Departamento de Matem\'atica\\
	Universidade Federal do Recôncavo da Bahia\\
	Federal University of Recôncavo da Bahia (URFB)- C\^ampus de Cruz das Almas \\
	44380-000, Cruz das Almas, BA, Brazil}
\email{carlosponcio@ufrb.edu.br}

\author[Edivaldo L. dos Santos]{Edivaldo L. dos Santos}
\address{Departamento de Matem\'atica\\
	Universidade Federal de S\^{a}o Carlos\\
	Federal University of S\~{a}o Carlos  (UFSCAR) - C\^ampus de S\~ao Carlos \\
	13565-905, S\~ao Carlos, SP, Brazil}
\email{edivaldo@ufscar.brr}

\author[Leandro V. Mauri]{Leandro V. Mauri}
\address{Departamento de Matem\'atica\\
	Instituto de Ci\^encias Matem\'aticas e de Computa\c c\~ao\\
	S\~ao Paulo University (USP)- C\^ampus de S\~ao Carlos \\
	13560-970, S\~ao Carlos, SP, Brazil}
\email{leandro.mauri@usp.br}

\author[ D. De Mattos]{Denise de Mattos}
\address{Departamento de Matem\'atica\\
	Instituto de Ci\^encias Matem\'aticas e de Computa\c c\~ao\\
	S\~ao Paulo University (USP)- C\^ampus de S\~ao Carlos \\
	13560-970, S\~ao Carlos, SP, Brazil}
\email{deniseml@icmc.usp.br}

\thanks{The first author was supported by FAPESP of Brazil Grant numbers 2017/21162-0 and 2016/05463-7.}

\thanks{The third author was supported by FAPESP of Brazil Grant number 2018/23928-2.}

\begin{abstract}In this paper, we present some results related to Bárány-Larman colored problem and The \v{Z}ivaljevi\'{c} and Vre\'{c}ica colored Tverberg problem. We give an alternative proof for the Bárány-Larman Conjecture for primes $-1$ and the optimal colored Tverberg theorem. Also, we prove generalizations of  the Colored Tverberg theorem of \v{Z}ivaljevi\'{c} and Vre\'{c}ica. \end{abstract}


\maketitle

\section{Introduction}

\vspace{0.2cm}

In 1959 Birch \cite{Bir59} formulated the following conjecture. 

\begin{conj} Any $(r-1) (d+1) +1$ points in $\mathbb{R}^{d}$ can be partitioned in $N$ subsets whose convex hulls have a common point.
\end{conj}

The Birch's conjecture was proved by Helge Tverberg (see \cite{Tve66}, \cite{Tve81}) and since then it is known as Tverberg theorem. 

\begin{teo}[\textbf{Tverberg theorem}] Let $d \ge 1$ and $r \ge 2$ be integers, and $N= (r-1) (d+1)$. For any affine map $f: \Delta_N \rightarrow \mathbb{R}^{d}$ there are $r$ pairwise disjoint faces $\sigma_1, \ldots , \sigma_r$ of $\Delta_N$ such that $f(\sigma_1) \cap \cdots \cap f(\sigma_r) \neq \emptyset.$

\end{teo}

The following supposition is a generalization of the Tverberg theorem to arbitrary continuous map.

\begin{conj}[\textbf{Topological Tverberg conjecture}] Let $d \ge 1$ and $r \ge 2$ be integers, and let us consider $N= (r-1) (d+1)$. For any continuous map $f: \Delta_N \rightarrow \mathbb{R}^{d}$ there are $r$ pairwise disjoint faces $\sigma_1, \ldots , \sigma_r$ of $\Delta_N$ such that $f(\sigma_1) \cap \cdots \cap f(\sigma_r) \neq \emptyset.$

\end{conj}

The Topological Tverberg conjecture was considered an important unsolved problem in topological combinatorics \cite{mat08}. In 1981 the conjecture was proved, when $r$ is a prime number, by Bárány, Shlosman and Szücs \cite{Bar81}. The result for a prime number $r$ was extend for a prime power $r$ by Özaydin (unpublished) \cite{Oza87}, Volovikov \cite{Vol96} and Sarkaria \cite{Sar00}. This result is known as Topological Tverberg theorem.

 \begin{teo}[\textbf{Topological Tverberg theorem}] Let $d \ge 1$, $r \ge 2$, and $N= (r-1) (d+1)$ be integers. If $r$ is a prime power, then for any continuous map $f: \Delta_N \rightarrow \mathbb{R}^{d}$ there are $r$ pairwise disjoint faces $\sigma_1, \ldots , \sigma_r$ of $\Delta_N$ such that $f(\sigma_1) \cap \cdots \cap f(\sigma_r) \neq \emptyset.$

\end{teo}

The set $\{\sigma_1,\ldots ,\sigma_r \}$ of disjoint faces of $\Delta_N$ whose images by the map $f: \Delta_N \rightarrow \mathbb{R}^{d}$ have nonempty intersection is called a Tverberg partition for $f$.

\vspace{0.1cm}

In 1992, Bárány and Larman \cite{Bar92} formulated the colored Tverberg problem\ as follows.

\begin{defi}[\textbf{Coloring}] Let $N \ge 1$ be a integer and let $V(\Delta_N)$ be the set of vertices of the simplex $\Delta_N$. A \textit{coloring} of vertices $V(\Delta_N)$ by $l$ colors is a partition $(C_1, \ldots , C_l)$ of $V(\Delta_N)$, that is, $V(\Delta_N) = C_1 \cup \cdots \cup C_l$ and $C_i \cap C_j = \emptyset$ for $1 \le i < j \le l$. The elements of the partition $(C_1, \ldots , C_l)$ are called \textit{color classes}.

\end{defi}

\begin{defi}[\textbf{Rainbow face}] Let $(C_1, \ldots , C_l)$ be the coloring of $V(\Delta_N)$ by $l$ colors. A face $\sigma$ of the  simplex $\Delta_N$ is a \textit{rainbow face} if $|\sigma \cap C_i| \le 1$, for all $1 \le i \le l$.

\end{defi}

\begin{prob}[\textbf{Bárány-Larman colored problem}] Let $d \ge 1$ and $r \ge 2$ be integers. Determine the smallest number $n= n(d,r)$ such that for every map $f: \Delta_{n-1} \rightarrow \mathbb{R}^{d}$, and every coloring $(C_1, \ldots , C_{d+1})$ of the vertex set $V(\Delta_{n-1})$ of the simplex $\Delta_{n-1}$ by $d+1$ colors, with each color of size at least $r$, there exist $r$ pairwise disjoint rainbow faces $\sigma_1, \ldots , \sigma_r$ of $\Delta_{n-1}$ satisfying \begin{eqnarray} \hspace{0.2cm} f(\sigma_1) \cap \cdots \cap f(\sigma_r) \neq \emptyset. \nonumber\end{eqnarray}

\end{prob}

\vspace{0.1cm}

Note that a trivial lower bound for $n(d,r)$ is $(d+1)r$. Then the following conjecture arises.

\vspace{0.1cm}

\begin{conj}[\textbf{Bárány-Larman Conjecture}] Let $d\ge 1$ and $r \ge 2$ be integers. Then $n(d,r)= (d+1)r$.

\end{conj}

\vspace{0.2cm}

In \cite{Bla15}, the following result was formulated and proved.

\vspace{0.2cm}

\begin{teo}[\textbf{The optimal colored Tverberg theorem}] Let $d \ge 1$ be an integer, let $r \ge 2$ be a prime number and $N \ge (d+1)(r-1)$. For every continuous map $f: \Delta_N \rightarrow \mathbb{R}^{d}$, and every coloring $(C_1, \ldots , C_{m})$ of the vertex set $V(\Delta_N)$ by $m$ colors, with each color of size at most $r-1$, there exist $r$ pairwise disjoint rainbow faces $\sigma_1, \ldots , \sigma_r$ of $\Delta_N$ satisfying: \begin{eqnarray} \hspace{0.2cm} f(\sigma_1) \cap \cdots \cap f(\sigma_r) \neq \emptyset. \nonumber\end{eqnarray}

\end{teo}

\vspace{0.2cm}

In the same paper \cite{Bla15}, it is proved that the optimal colored Tverberg theorem implies the Bárány-Larman Conjecture, for the case when $r+1$ is a prime, which is called Bárány-Larman Conjecture for primes --1.

\vspace{0.1cm}

In this paper, we give a shorter alternative proof, without using spectral sequences, for the fact that $\mbox{Index}_{\mathbb{Z}_p}S(W_p)=\mbox{Index}_{\mathbb{Z}_p}\Delta _{p,p-1}=H^{\geq p-1}( B\mathbb{Z}_{p})$ (where $p$ is a prime number). This is the key result to prove the optimal colored Tverberg theorem and the Bárány-Larman Conjecture for primes -$1$

\vspace{0.1cm}

A modified colored Tverberg problem was presented by \v{Z}ivaljevi\'{c} and Vre\'{c}ica in the paper \cite{Ziv92}.

\begin{prob}[\textbf{The \v{Z}ivaljevi\'{c} and Vre\'{c}ica colored Tverberg problem}] Let $d \ge 1$ and $r \ge 2$ be integers. Determine the smallest number $t=t(d,r)$ such that for every affine (or continuous) map \linebreak $f: \Delta \rightarrow \mathbb{R}^{d}$, and every coloring $(C_1, \ldots , C_{d+1})$ of the the vertex set $V(\Delta)$ by $d+1$ colors, with each color of size at least $t$, there exist $r$ pairwise disjoint rainbow faces $\sigma_1, \ldots , \sigma_r$ of $\Delta_{n-1}$ satisfying \begin{eqnarray}  f(\sigma_1) \cap \cdots \cap f(\sigma_r) \neq \emptyset.\nonumber\end{eqnarray}

\end{prob}

For $r \ge 2$ be a prime power, \v{Z}ivaljevi\'{c} and Vre\'{c}ica proved that $t(d,r) \le 2r-1$. This result was known as colored Tverberg theorem of \v{Z}ivaljevi\'{c} and Vre\'{c}ica.

\begin{teo}[\textbf{Colored Tverberg theorem of \v{Z}ivaljevi\'{c} and Vre\'{c}ica} \cite{Ziv92}] Let $d \ge 1$ be an integer, and let $r \ge 2$ be a prime power. For every continuous map $f: \Delta \rightarrow \mathbb{R}^{d}$, and every coloring $(C_1, \ldots , C_{d+1})$ of the the vertex set $V(\Delta)$ by $d+1$ colors, with each color of size at least $2r-1$, there exist $r$ pairwise disjoint rainbow faces $\sigma_1, \ldots , \sigma_r$ of $\Delta$ satisfying: \begin{eqnarray} f(\sigma_1) \cap \cdots \cap f(\sigma_r) \neq \emptyset.\nonumber\end{eqnarray}

\end{teo}

\vspace{0.2cm}

We prove in this paper that the Colored Tverberg theorem of \v{Z}ivaljevi\'{c} and Vre\'{c}ica holds when $| C_1 | \ge 2r-1$ and $|C_i| \ge 2r-4$, for all $i=2,\ldots , d+1$. We also prove that this theorem still holds with a few more variations in the cardinalities of the color classes.

\vspace{1cm}

\section{Preliminares}

In this section we present the preliminary results.

\vspace{0.2cm}

\subsection{Volovikov index}

In this section, we present the Volovikov index, which is an algebraic topology tool used in the proofs of the theorems in Sections 4 and 5. For more details on the Volovikov index, see \cite{Vol00}.

\vspace{0.2cm}

Initially, we define the \textit{Borel construction}.

\vspace{0.2cm}

\begin{defi}[\textbf{Borel construction}] Let $G$ be a compact Lie group and let $X$ be a Hausdorff paracompact  $G$-space on which $G$ acts freely. Then $X \rightarrow X / G$ is a principal $G$-bundle (see \cite{Bre72}) and one can take:  \begin{eqnarray} h: X / G \rightarrow BG\nonumber\end{eqnarray} a classifying map for the $G$-bundle $X \rightarrow X / G$.

\vspace{0.1cm}

Let us consider the product $EG \times X$ with diagonal $G$-action $g(e,x) =(ge, gx)$, for every $g \in G$ and $(e,x) \in EG \times X$. Let $EG \times_G X = (EG \times X) / G$ be its orbit space. The first projection $EG \times X \rightarrow EG$ induces a map: \begin{eqnarray} p_X: EG \times_G X \rightarrow (EG)/G = BG, \nonumber \\ (e,x) G \mapsto e G  \hspace{1.9cm} \nonumber \end{eqnarray} which is a fibration with fiber $X$ and base space $BG$ being the classifying space of $G$.

\vspace{0.1cm}

This is the \textit{Borel construction}. It associates to each Hausdorff paracompact $G$-space $X$, a $G$-space $EG \times_G X$, which is denoted by $X_G$, over $BG$. Also, it associates to each $G$-map $f: X \rightarrow Y$, a fiber preserving map $\overline{ \mbox{Id}_{EG} \times f}:EG \times_G X \rightarrow EG \times_G Y$.

\end{defi}

Let us recall the Leray-Serre theorem for fibrations.

\begin{teo}[\textbf{The cohomology Leray-Serre Spectral sequence (Theorem $5.2$ \cite{McC01})}]
Let $R$ be a commutative ring with unit. Given a fibration $F \hookrightarrow E \stackrel{p}{\rightarrow} B$, where $B$ is a path-wise connected space, there is a first quadrant spectral sequence of algebras $\{ E_r^{\ast, \ast}, d_r \}$, with: \begin{eqnarray} E_{2}^{p,q} \cong H^{p} ( B ; \mathcal{H}^{q} (F;R)), \nonumber \end{eqnarray} the cohomology of $B$, with local coefficients in the cohomology of $F$, the fiber of $p$, and converging to $H^{\ast}(E;R)$ as an algebra. Furthermore, this spectral sequence is natural with the respect to fiber-preserving maps of fibrations.

\end{teo}

Let us now recall one of the numerical index defined by Volovikov in \cite{Vol00}. This is a function on $G$-spaces whose value is either a positive integer or $\infty$. For our purposes, it is sufficient to consider that $G$ is a $p$-torus ($p$ a prime number), that is, $G= (\mathbb{Z}_p)^{n}$, for $n \ge 1$, and $\mathbb{K} = \mathbb{Z}_p$ is a coefficient field in \v{C}ech cohomology .

\vspace{0.1cm}

\begin{defi}[\textbf{Volovikov index}] Let $G= (\mathbb{Z}_p)^{n}$ ($p$ a prime number) be a compact Lie group, for $n \ge 1$, and let $X$ be a Hausdorff paracompact  $G$-space, on which $G$ acts freely. The definition of the \textit{Volovikov index of $X$}, denoted by $i(X)$, uses the spectral sequence of the bundle $p_X: X_G \rightarrow BG$, with fibre $X$ (the Borel construction), given in Theorem $2.2$. This spectral sequence converges to the equivariant cohomology $H^{\ast} (X_G; \mathbb{Z}_p)$. Let $\Lambda^{\ast}$ be the equivariant cohomology algebra of a point $H^{\ast} (\mbox{pt}_G; \mathbb{Z}_p)= H^{\ast} (BG; \mathbb{Z}_p)$. Suppose that $X$ is path connected. Then $E_{2}^{\ast,0} = \Lambda^{\ast}$. Assume that $E_{2}^{\ast,0}= \cdots = E_{s}^{\ast,0} \neq E_{s+1}^{\ast,0}$. Then, by definition, $i(X)=s$. If $E_{2}^{\ast,0}= \cdots = E_{\infty}^{\ast,0}$ then, by definition, $i(X)= \infty$.

\end{defi}

We state some properties of the Volovikov index (see \cite{Vol00}).

\begin{prop} Let $X$, $Y$ and $Z$ be a Hausdorff paracompact  $G$-space where the $G$ acts freely. Then

\vspace{0.2cm}
\begin{itemize}
\item[(i)]  If there is a $G$-equivariant map $X \rightarrow Y$, then $i(X) \le i(Y)$.

\item[(ii)]  If $\tilde{H}^{j}(X; \mathbb{Z}_p)=0$, for all $j<n$, then $i(X) \ge n+1$.

\item[(iii)]  If $H^{j}(Z; \mathbb{Z}_p)=0$, for all $j \ge n$ and if $i(Z)< \infty$, then $i(Z) \le n$.

\item[(iv)]   If $X$ is a compact or finite dimensional space such that $H^{\ast}(X; \mathbb{Z}_p) = H^{\ast}(S^{n}; \mathbb{Z}_p) $, and if $G$ acts on $X$ without fixed points, then $i(X) =n+1$.

\end{itemize}

\end{prop}

\vspace{0.3cm}

\subsection{Fadell-Husseini index}

Next, we introduce the Fadell-Husseini index, and some of its properties.

\vspace{0.2cm}

\begin{defi}[\textbf{Fadell-Husseini index}] Let $G$ be a finite group and let $R$ be a commutative ring with unit. For a $G$-space $X$, the \textbf{Fadell-Husseini index} of $X$, denoted by $\mbox{Index}_G (X;R)$, is defined as follow: \begin{eqnarray} \mbox{Index}_G (X;R) = \mbox{Ker } ( q_{X}^{\ast}: H^{\ast} (BG;R) \rightarrow H^{\ast} (EG \times_G X;R) ), \nonumber \end{eqnarray}
where $q_X: X / G \rightarrow BG$ is a classifying map for the $G$-bundle $X \rightarrow X / G $.

\end{defi}

\vspace{0.2cm}

\begin{prop}[monotonicity property] If $f: X \rightarrow Y$ is a $G$-equivariant map, then \begin{eqnarray} \mbox{Index}_G (Y;R) \subset \mbox{Index}_G (X;R). \nonumber \end{eqnarray}

\end{prop}

\vspace{0.2cm}

\begin{defi} $W_r= \{ (t_1, \cdots, t_r) \in \mathbb{R}^{r} \hspace{0.1cm}; \hspace{0.1cm} \sum_{i=1}^{r} t_i =0\}.$
\end{defi}

\vspace{0.2cm}

\begin{teo}[\cite{Bla17}, pg. 305]

Let $r \ge 1$ be an integer, $p \ge 2$ a prime number and $N=(r-1) (d+1)$. Then: \begin{eqnarray} \mbox{Index}_{\mathbb{Z}_p} \left(S (W_{r}^{\oplus (d+1)}; \mathbb{Z}_p\right) = \langle b \rangle \mbox{, where } b \in H^{N}(B \mathbb{Z}_p). \nonumber \end{eqnarray}

\end{teo}

\vspace{0.3cm}

\subsection{Chessboard complex and connectedness}

In this section, we introduce the chessboard complex, which is widely used in Colored Tverberg theorem proofs. Also, we introduce the concept of connectedness and we relate it to the Volovikov index. For more details on these topics, see \cite{Bla17} and \cite{mat08}.

\vspace{0.3cm}

\begin{defi}[\textbf{Chessboard complex}] The \textit{$m \times n$ chessboard complex} $\Delta_{m,n}$ is the simplicial complex whose vertex set is $[m] \times [n]$\footnote{$[n]=\{1, \cdots,n\}.$} . The simplexes of $\Delta_{m,n}$ are the subsets $\{(i_0,j_0), \ldots, (i_k,j_k) \} \subset [m] \times [n]$, where $[n]=\{1, \ldots, n\}$, $i_s \neq i_{s'}$ ($1 \le s < s'\le k$), and $j_t \neq j_{t'}$ ($1 \le t < t'\le k$).

\end{defi}

\vspace{0.1cm}

\begin{defi} Let $n \ge -1$ be an integer. A topological space $X$ is \textit{$n$-connected} if any continuous map $f: S^{k} \rightarrow X$, where $-1 \le k \le n$, can be continuously extended to a continuous map $g: B^{k+1} \rightarrow X$, that is $g|_{\partial B^{k+1} = S^k} = f$ (here $B^{k+1}$ denotes a $(k+1)$-dimensional closed ball whose boundary is the sphere $S^{k}$). A topological space is $(-1)$-connected if it is non-empty. If the space $X$ is $n$-connected, but it is not $(n+1)$-connected, we write $\mbox{conn} (X)=n$.

\end{defi}

\vspace{0.1cm}

\begin{teo}[pg. 332, \cite{Bla17}] Let $X$ and $Y$ be topological spaces. Then: \begin{eqnarray}\mbox{conn} (X \ast Y) \ge  \mbox{conn} (X) + \mbox{conn} (Y) +2. \hspace{1cm} \nonumber \end{eqnarray} \end{teo}

\vspace{0.1cm}

\begin{teo}[ 4.4.1 Theorem \cite{mat08}] Let $X$ be a nonempty topological space and let $k \ge 1$. Then $X$ is $k$-connected, if only if, it is simply connected (i.e., the fundamental group $\pi_1(X)$ is trivial) and $\tilde{H}_i(X)=0$, for all $i=0,1, \ldots,k$.

\end{teo}

\vspace{0.1cm}

\begin{teo} Let be a $X$ topological space. Then \begin{eqnarray} i (X ) \ge  \mbox{conn} (X) +2. \nonumber \end{eqnarray}
\end{teo}

\vspace{0.1cm}

\begin{teo} Let $m, n \ge 1$ be integers. Then \footnote{$\lfloor{x}\rfloor = \mbox{max } \{ n \in \mathbb{N} \mbox{; } n \le x \}.$ } \begin{eqnarray} \mbox{conn} (\Delta_{m,n}) = \mbox{min} \{m,n, \lfloor{\frac{m+n+1}{3}}\rfloor\}-2. \nonumber \end{eqnarray}

\end{teo}

\vspace{0.1cm}

\begin{cor} Let $r \ge 2$ be an integer. Then $\mbox{conn} (\Delta_{2r-1,r}) = r-2$.

\end{cor}

\vspace{1cm}

\section{Alternative proof for the optimal colored Tverberg theorem}

In this section, we provide an alternative proof for the Theorem $3.3$ (\cite{Bla17} Theorem $6.6$), which is the key result used in the proof of B\'{a}r\'{a}ny-Larman Conjecture, for primes --1, and the optimal colored Tverberg theorem. Our proof is shorter than that shown in \cite{Bla17} and it does not use spectral sequences.

\vspace{0.2cm}

Next, we present the two results that are used in the proof of Theorem $3.3$.

\vspace{0.2cm}

\begin{lemma}[\textbf{Lemma 6.5 \cite{Bla17}}]
Let $p$ be an odd prime. There exists a $\mathbb{Z}_p$-equivariant map $f : \Delta_{p-1,p}\rightarrow S(W_p)$ such that the induced map in cohomology $f^* :  H^{p-2}(S(W_p);\mathbb{Z}_p)\rightarrow H^{p-2}(\Delta_{p-1,p};\mathbb{Z}_p)$ is an isomorphism.
\end{lemma}

\vspace{0.2cm}

\begin{lemma}[\textbf{Lemma 3 \cite{Gon04}}]

Let $X$ be a Hausdorff, connected and locally pathwise connected space and let $G$ be a finite group acting freely on $X$. Then for any $i \ge 0$ and a commutative ring $R$ with a unit, there is a transfer homomorphism: 

$$\tau_X: H^i(X;R)\longrightarrow H^i(X/ G;R)$$
with the following properties:

\begin{itemize}
\item[(i)]If $X$ is a $k$-dimensional $G-CW$-complex (not necessarily finite), then \begin{eqnarray} \tau_X: H^k(X;R)\longrightarrow H^k(X/ G;R) \mbox{ is surjective;} \nonumber \end{eqnarray}

\item[(ii)] If $Y$ is another space that satisfies the hypothesis of the lemma, $h:X\rightarrow Y$ is an equivariant map and $\overline{h} : X/G \rightarrow Y /G$ is the map induced by $h$, then $\tau_X \circ h^*=\overline{h} ^*\circ\tau_Y$.
\end{itemize}

\centerline{ \xymatrix{ H^{i}(Y) \ar[d]_{\tau_Y} \ar[r]^{h^{\ast}} & H^{i} (X) \ar[d]^{\tau_X} \\ H^{i} (Y/G) \ar[r]^{\overline{h^{\ast}}} & H^{i} (X/G) }  }

\end{lemma}

\vspace{0.2cm}

Next, we prove the main result in this section.

\vspace{0.2cm}

\begin{teo}[\textbf{Theorem 6.6 \cite{Bla17}}]

Let $p$ be an prime number. Then: \begin{eqnarray}\mbox{Index}_{\mathbb{Z}_p}S(W_p)=\mbox{Index}_{\mathbb{Z}_p}\Delta _{p,p-1}=H^{\geq p-1}( B\mathbb{Z}_{p}). \nonumber \end{eqnarray}

\end{teo}

\textit{Proof.} 
Let $f:\Delta_{p,p-1} \rightarrow S(W_p)$ be the $\mathbb{Z}_p$-equivariant map of the Lemma $3.1$ and let $\overline{f}: \Delta_{p,p-1}/ \mathbb{Z}_p \rightarrow S(W_p)/ \mathbb{Z}_p$  be the induced map in the orbit space.

\vspace{0.2cm}

If $q': S(W_p)/ \mathbb{Z}_p \rightarrow B \mathbb{Z}_p$ is a classifying map for the $\mathbb{Z}_p$- principal bundle $S(W_p) \rightarrow S(W_p) / \mathbb{Z}_p$, then $q' \circ \overline{f} : \Delta_{p,p-1} / \mathbb{Z}_p \rightarrow B \mathbb{Z}_p$ is a classifying map for the $\mathbb{Z}_p$- principal bundle $\Delta_{p,p-1} \rightarrow \Delta_{p,p-1} / \mathbb{Z}_p$.

\vspace{0.2cm}

Let us consider the transfer homomorphisms as in the Lemma $3.2$:

\begin{eqnarray}  \tau_{S(W_p)}: H^{p-2} (S(W_p)) \rightarrow H^{p-2} (S(W_p) / \mathbb{Z}_p) \nonumber  \\ \tau_{\Delta_{p,p-1}}: H^{p-2} (\Delta_{p,p-1}) \rightarrow H^{p-2} (\Delta_{p,p-1} / \mathbb{Z}_p)  \nonumber \end{eqnarray}

\vspace{0.2cm}

Since $f$ is a $\mathbb{Z}_p$-equivariant map, by Lemma $3.2 (ii)$, we have the following commutative diagram,

\centerline{ \xymatrix{H^{p-2} (S(W_p)) \ar[r]^{f^{\ast}} \ar[d]_{\tau_{S(W_p)}} & H^{p-2} (\Delta_{p,p-1}) \ar[d]^{\tau_{\Delta_{p,p-1}}} \\ H^{p-2} (S(W_p)/ \mathbb{Z}_p) \ar[r]^{\overline{f}^{\ast}}  & H^{p-2} (\Delta_{p,p-1}/ \mathbb{Z}_p) } }

\noindent that is, \begin{eqnarray} \tau_{\Delta_{p,p-1}} \circ f^{\ast} = \overline{f}^{\ast} \circ \tau_{S(W_p)}.\end{eqnarray}

\vspace{0.2cm}

It follows from Lemma $3.1$ that $f$ is an isomorphism. By Lemma $3.2 (i)$, the tranfer homomorphism $\tau_{p,p-1}$ is surjective, thus the composition $\tau_{\Delta_{p,p-1}} \circ f^{\ast}$ is not the null homomorphism, namely \begin{eqnarray} \tau_{\Delta_{p,p-1}} \circ f^{\ast} \neq 0. \end{eqnarray} 

\vspace{0.2cm}

Hence, from $(1)$ and $(2)$, we have $0 \neq \tau_{\Delta_{p,p-1}} \circ f^{\ast} = \overline{f}^{\ast} \circ \tau_{S(W_p)}$. Then $\overline{f}^{\ast}$ is not trivial, that is, \begin{eqnarray} \overline{f}^{\ast} \neq 0.\end{eqnarray}

\vspace{0.2cm}

Let us consider the following commutative diagram

\centerline{ \xymatrix{H^{p-2} (S(W_p)) \ar[rr]^{f^{\ast}} \ar[d]_{\tau_{S(W_p)}} &  &H^{p-2} (\Delta_{p,p-1}) \ar[d]^{\tau_{\Delta_{p,p-1}}} \\ H^{p-2} (S(W_p)/ \mathbb{Z}_p) \ar[rr]^{\overline{f}^{\ast}}  &  &H^{p-2} (\Delta_{p,p-1}/ \mathbb{Z}_p)  \\ & H^{p-2} (B \mathbb{Z}_p) \ar[lu]^{q'^{\ast}} \ar[ru]_{q^{\ast}}}  }

\vspace{0.1cm}

\noindent it means that $q^{\ast}= \overline{f}^{\ast} \circ q'^{\ast}$, where $q$ and $q'$ are classifying maps.
 
\vspace{0.2cm}

By Theorem $2.7$, we have that $\mbox{Index}_{\mathbb{Z}_p} S(W_p) = H^{ \ge p-1} ( B \mathbb{Z}_p)$ and therefore \begin{eqnarray}
\mbox{the homomorphism $q'^{\ast}$ \mbox{at level} $p-2$ \mbox{ is an isomorphism.}}\end{eqnarray} 

\vspace{0.2cm}

As a consequence from $(3)$, $(4)$ and by the last commutative diagram, we have that \begin{eqnarray}
\mbox{the homomorphism $q^{\ast}$ \mbox{at level} $p-2$ \mbox{ is an isomorphism.}} \nonumber\end{eqnarray}

\vspace{0.2cm}

Also, it is a known fact that\begin{eqnarray}
H^{\ast} (B \mathbb{Z}_p ) \cong \mathbb{Z}_p [y] \otimes_{\mathbb{Z}_p} \Lambda (x), \mbox{where } y \in H^{2}(B \mathbb{Z}_p),\hspace{0.1cm} x \in H^{1} (B \mathbb{Z}_p) \mbox{ and } x^{2}=0, \nonumber \end{eqnarray} 

\vspace{0.2cm}

where $\Lambda ( \cdot )$ denotes the exterior algebra.

\vspace{0.2cm}

The generator $\mu \in H^{i} (B \mathbb{Z}_p)$ is of the form 

\begin{eqnarray}\mu=\left\{\begin{array}{rc}
x y^{\frac{(i-2)}{2}},&\mbox{if}\quad i \mbox{ is odd},\\
y^{\frac{1}{2}}, &\mbox{if}\quad i \mbox{ is even}.
\end{array}\right. \nonumber
\end{eqnarray}

\vspace{0.2cm}

Then, we have that: \begin{eqnarray} q^{\ast}(x y^{\frac{p-3}{2}}) \neq 0.\end{eqnarray}

\vspace{0.2cm}

Take the $\mathbb{Z}_p$-equivariant map $f:\Delta_{p-1,p} \rightarrow S(W_p))$ of the Lemma $3.1$. By Theorem $2.7$, we have: \begin{eqnarray} \mbox{Index}_{\mathbb{Z}_p} S(W_p) = H^{ \ge p-1} ( B \mathbb{Z}_p) \subset \mbox{Index}_{\mathbb{Z}_p} (\Delta_{p-1,p}).  \end{eqnarray}

\vspace{0.2cm}

Let $q^{\ast}: H^{i} (B \mathbb{Z}_p) \rightarrow H^{i} (\Delta_{p-1,p}/ \mathbb{Z}_p)$ be the classifying map, where $i< p-2$. Next, we prove that $Ker (q^{\ast})=0$, for $i< p-2$. 

\vspace{0.2cm}

It is sufficient to prove that $q^{\ast}(\mu) \neq 0$, where $\mu$ is the generator of $H^{i}(B \mathbb{Z}_p)$.

\vspace{0.2cm}

Let us suppose that $i$ is an odd number and $q^{\ast}(\mu)= q^{\ast}(x y^{\frac{(i-1)}{2}})=0$. Then: \begin{eqnarray} q^{\ast} \left(x y ^{\frac{(p-3)}{2}}\right) = q^{\ast} \left( x y ^{\frac{(i-1)}{2}} y^{\frac{(p-i-2)}{2}}\right)= \\ q^{\ast} \left( x y ^{\frac{(i-1)}{2}}\right) \hspace{0.05cm} q^{\ast}\left( y^{\frac{(p-i-2)}{2}}\right)=0. \nonumber \end{eqnarray}

But, $(7)$ contradicts $(5)$.

\vspace{0.2cm}

Now, let us suppose that $i$ is an even number and $q^{\ast}(\mu)=q^{\ast} (y^{\frac{i}{2}})=0$. Then: \begin{eqnarray} q^{\ast} \left(x y ^{\frac{(p-3)}{2}}\right) = q^{\ast} \left( x y ^{\frac{(p-i-3)}{2}} y^{\frac{i}{2}}\right)= \\ q^{\ast} \left( x y ^{\frac{(p-i-3)}{2}}\right) \hspace{0.05cm} q^{\ast}\left( y^{\frac{i}{2}}\right)=0. \nonumber \end{eqnarray}

But, $(8)$ contradicts $(5)$.

\vspace{0.2cm}

Therefore,
\begin{eqnarray} \mbox{Index}_{\mathbb{Z}_p} (\Delta_{p-1,p}) = Ker(q^{\ast}) \subset H^{ \ge p-1} ( B \mathbb{Z}_p).  \end{eqnarray}

\vspace{0.2cm}

Consequently, from $(6)$ and $(9)$, we have that: \begin{eqnarray}\mbox{Index}_{\mathbb{Z}_p} (\Delta_{p-1,p}) = H^{ \ge p-1} ( B \mathbb{Z}_p).\nonumber \end{eqnarray}.

\hspace{14.3cm} $\square$

\vspace{0.2cm}

Let us note that Theorem $3.3$ is essential in the proof of the following results (for details, see \cite{Bla17}).

\vspace{0.2cm}

\begin{teo}[\textbf{Bárány-Larman Conjecture for primes $-1$}] Let $d\ge 1$ and $r \ge 2$ be integers such that $p=r+1$ is a prime number. Then $n(d,r)= (d+1)r$.

\end{teo}

\vspace{0.2cm}

\begin{teo}[\textbf{The optimal colored Tverberg theorem}] Let $d \ge 1$ be an integer, let $r \ge 2$ be a prime and $N \ge (d+1)(r-1)$. For every continuous map $f: \Delta_N \rightarrow \mathbb{R}^{d}$, and every coloring $(C_1, \ldots , C_{m})$ of the vertex set $V(\Delta_N)$ by $m$ colors, with each color of size at most $r-1$, there exist $r$ pairwise disjoint rainbow faces $\sigma_1, \ldots , \sigma_r$ of $\Delta$ satisfying \begin{eqnarray}  f(\sigma_1) \cap \cdots \cap f(\sigma_r) \neq \emptyset.\nonumber\end{eqnarray}

\end{teo}

\vspace{0.5cm}

\section{Generalization of the Colored Tverberg theorem of \v{Z}ivaljevi\'{c} and Vre\'{c}ica}

\vspace{0.2cm}

In this section, we prove a generalization of the Colored Tverberg theorem of \v{Z}ivaljevi\'{c} and Vre\'{c}ica (\cite{Ziv92}), in which is assumed only one color class with cardinality at least $2r-1$ and the other $d$ color classes with  cardinality cardinality $2r-4$.

\vspace{0.2cm}

First, we state the following result which is used in the proof of our main theorem (Theorem $4.2$).

\vspace{0.2cm}

\begin{teo}[\textbf{Corollary $4.8$ \cite{Bla17}}] Let $d \ge 1$ and $m \ge 1$ be integers, and let $r =p^{n} \ge 2$ be a prime power. Let $(C_1, \ldots , C_m)$ be a coloring of $\Delta$ by $m$ colors. If there is no $\mathfrak S_{r}$-equivariant map: \begin{eqnarray} \Delta_{|C_1|,r} \ast \cdots \ast \Delta_{|C_m|,r} \longrightarrow S \left(W_{r}^{\oplus (d+1)}\right)\nonumber \end{eqnarray} then for every continuous map $f: \Delta \rightarrow \mathbb{R}^{d}$ there exist $r$ pairwise disjoint rainbow faces $\sigma_1, \ldots, \sigma_r$ of $\Delta$ satifying \begin{eqnarray}f(\sigma_1) \cap \cdots \cap f(\sigma_r) \neq \emptyset.\nonumber\end{eqnarray}

\end{teo}

\vspace{0.2cm}

Now, we can formulate and prove the main result of this section.

\vspace{0.3cm}

\begin{teo}[\textbf{Generalization of the Colored Tverberg theorem of \v{Z}ivaljevi\'{c} and Vre\'{c}ica}] Let $d \ge 1$ be an integer, and let $r \ge 2$ be a prime power. For every continuous map $f: \Delta \rightarrow \mathbb{R}^{d}$, and every coloring $(C_1, \ldots , C_{d+1})$ of the the vertex set $V(\Delta)$ by $d+1$ colors, with $|C_1| \ge 2r-1$ and $|C_i| \ge 2r-4$, for all $i =2, \cdots, d+1$, there exist $r$ pairwise disjoint rainbow faces $\sigma_1, \ldots , \sigma_r$ of $\Delta$ satisfying: \begin{eqnarray} f(\sigma_1) \cap \cdots \cap f(\sigma_r) \neq \emptyset.\nonumber\end{eqnarray}

\end{teo}

\vspace{0.2cm}

\textit{Proof.} We split the proof into two cases. 

\vspace{0.2cm}

\textbf{1º Case: $r \ge 3:$} Without loss of generality, we can assume that $|C_1|=2r-1$ and $|C_2| = \cdots = |C_m| = 2r-4$. Then, by Theorem $4.1$, it is sufficient to show that there is no $\mathfrak S_{r}$-equivariant map: \begin{eqnarray}\Delta_{2r-1,r} \ast  (\Delta_{2r-4,r})^{\ast d} \longrightarrow S \left(W_{r}^{\oplus (d+1)}\right). \nonumber \end{eqnarray}

\vspace{0.2cm}

Consider the regular embedding \begin{eqnarray} \varphi: (\mathbb{Z}_p)^{n} \longrightarrow \mbox{Sym } ((\mathbb{Z}_p)^{n}) \cong \mathfrak S_{r} \hspace{2.7cm} \nonumber \\ g \longmapsto L_g : (\mathbb{Z}_p)^{n} \rightarrow (\mathbb{Z}_p)^{n}, \hspace{0.1cm} L_g (x) =g+x.  \nonumber\end{eqnarray}

Therefore, we have a subgroup $(\mathbb{Z}_p)^{n} \cong \mbox{Im} (\varphi) \le \mathfrak S_{r}$. Then, to prove the non-existence of an $\mathfrak S_{r}$-equivariant map it suffices to prove the non-existence of an  $(\mathbb{Z}_p)^{n}$-equivariant map: \begin{eqnarray} \Delta_{2r-1,r} \ast (\Delta_{2r-4,r})^{\ast d} \longrightarrow S \left(W_{r}^{\oplus (d+1)}\right).\nonumber \end{eqnarray}

\vspace{0.2cm}

Let us suppose that there is an $(\mathbb{Z}_p)^{n}$-equivariant map. Then, by Proposition $2.4 (i)$, we have: \begin{eqnarray} i\left(\Delta_{2r-1,r} \ast (\Delta_{2r-4,r})^{\ast d} \right) \le i\left( S \left(W_{r}^{\oplus (d+1)}\right) \right). \nonumber\end{eqnarray}

\vspace{0.2cm}

Moreover, by Proposition $2.4(iv)$ we have: \begin{eqnarray}i\left( S \left(W_{r}^{\oplus (d+1)}\right) \right) = i\left( S^{(r-1) (d+1)-1}\right) = (r-1) (d+1). \nonumber \end{eqnarray}

\vspace{0.2cm}

We can calculate $\mbox{conn} (\Delta_{2r-1,r})$ and $\mbox{conn} (\Delta_{2r-4,r})$ using Theorem $2.13$, as follows.

 \begin{eqnarray} \mbox{conn} (\Delta_{2r-1,r})= \mbox{min} \left \{  2r-1, r, \left \lfloor \frac{(2r-1)+r +1}{3} \right \rfloor \right \}-2 = r-2 \mbox{  , and} \nonumber \end{eqnarray}

\begin{eqnarray}\mbox{conn} (\Delta_{2r-4,r})= \mbox{min} \left \{  2r-4, r, \left \lfloor \frac{(2r-4)+r +1}{3} \right \rfloor \right \}-2 = \nonumber \\ \mbox{min} \left \{  2r-4, r, r-1 \right \} \stackrel{(r \ge 3)}{=} (r-1)-2= r-3.\hspace{1cm} \nonumber \end{eqnarray}

\vspace{0.2cm}

It follows from Theorem $2.10$ that: \begin{eqnarray} \mbox{conn} \left(\Delta_{2r-1,r} \ast (\Delta_{2r-4,r})^{\ast d}\right) \ge  (r-2) + d (r-3) = (d+1) (r-1) -1.  \nonumber \end{eqnarray}

\vspace{0.2cm}

Hence, by Theorem $2.12$ we have: \begin{eqnarray} i \left(\Delta_{2r-1,r} \ast  (\Delta_{2r-4,r})^{\ast d}\right) \ge [(d+1)(r-1)-1]+2 = \nonumber \\ (d+1)(r-1)+1 > (r-1) (d+1). \hspace{2cm} \nonumber \end{eqnarray}

\vspace{0.3cm}

Then, $i\left(\Delta_{2r-1,r} \ast  (\Delta_{2r-4,r})^{\ast d}\right) > (r-1) (d+1) =i \left(S\left(W_{r}^{\oplus(d+1)}\right)\right)$.But it contradicts \begin{eqnarray}i\left(\Delta_{2r-1,r} \ast (\Delta_{2r-4,r})^{\ast d}\right) \le i \left(S\left(W_{r}^{\oplus(d+1)}\right)\right). \nonumber \end{eqnarray} 
Consequently, there is no $(\mathbb{Z}_p)^n$-equivariant map \begin{eqnarray}\Delta_{2r-1,r} \ast (\Delta_{2r-4,r})^{\ast d} \longrightarrow S \left(W_{r}^{\oplus (d+1)}\right). \nonumber \end{eqnarray}

\vspace{0.2cm}

\textbf{2º Case: $r =2:$} Note that from Theorem $2.13$, we have $\mbox{conn} (\Delta_{3,2})=0$ and $\mbox{conn} (\Delta_{1,2})=-1$. Thus, analogously to the first case, it follows that $i\left(\Delta_{3,2} \ast (\Delta_{1,2})^{\ast d}\right) > (d+1) (r-1)$.Therefore, the result holds for this case too.

\hspace{14.3cm} $\square$

\vspace{0.3cm}

\section{Flexibilization of Theorem $4.2$}

\vspace{0.2cm}

In this section, we present a flexibilization of Theorem $4.2$, in which we consider more variations in the minimum cardinalities for the $d+1$ color classes.

\vspace{0.3cm}

\begin{teo}[\textbf{Flexibilization of Theorem $4.2$}]  Let $d \ge 1$ be an integer, and let $r \ge 2$ be a prime power. For every continuous map $f: \Delta \rightarrow \mathbb{R}^{d}$, and every coloring $(C_1, \ldots , C_{d+1})$ of the the vertex set $V(\Delta)$ by $d+1$ colors, with  $|C_i| \ge 2r-1- 3x_i$  (where $x_i \ge 0$ is an integer), for all $i =1, \ldots, d+1$ such that $\sum_{i=1}^{d+1} x_i \le d$ and $2 x_i +1 \le r$, $ \forall i =1, \ldots, d+1$, there exist $r$ pairwise disjoint rainbow faces $\sigma_1, \ldots , \sigma_r$ of $\Delta$ satisfying \begin{eqnarray} f(\sigma_1) \cap \cdots \cap f(\sigma_r) \neq \emptyset.\nonumber\end{eqnarray}

\end{teo}

\vspace{0.2cm}

\textit{Proof.}  Without loss of generality, we can assume that $|C_i| =2r-1 -3x_i$, for all $i \in [d+1]$. Then, by Theorem $4.1$ it is sufficient to show that there is no $\mathfrak S_{r}$-equivariant map: \begin{eqnarray}\Delta_{2r-1-3x_1,r} \ast \cdots \ast \Delta_{2r-1-3x_i,r} \ast \cdots \ast \Delta_{2r-1-3x_{d+1},r} \longrightarrow S \left(W_{r}^{\oplus (d+1)}\right). \nonumber \end{eqnarray}

\vspace{0.2cm}

Let us consider the regular embedding \begin{eqnarray} \varphi: (\mathbb{Z}_p)^{n} \longrightarrow \mbox{Sym } ((\mathbb{Z}_p)^{n}) \cong \mathfrak S_{r} \hspace{2.7cm} \nonumber \\ g \longmapsto L_g : (\mathbb{Z}_p)^{n} \rightarrow (\mathbb{Z}_p)^{n}, \hspace{0.1cm} L_g (x) =g+x.  \nonumber\end{eqnarray}

Thus, we have a subgroup $(\mathbb{Z}_p)^{n} \cong \mbox{Im} (\varphi) \le \mathfrak S_{r}$. Hence, to prove the non-existence of an $\mathfrak S_{r}$-equivariant map it suffices to prove the non-existence of an  $(\mathbb{Z}_p)^{n}$-equivariant map: \begin{eqnarray}\Delta_{2r-1-3x_1,r} \ast \cdots \ast \Delta_{2r-1-3x_i,r} \ast \cdots \ast \Delta_{2r-1-3x_{d+1},r} \longrightarrow S \left(W_{r}^{\oplus (d+1)}\right). \nonumber \end{eqnarray}

\vspace{0.2cm}

Suppose that there is an $(\mathbb{Z}_p)^{n}$-equivariant map. Then, by Proposition $2.4 (i)$, we have: \begin{eqnarray}i\left(\Delta_{2r-1-3x_1,r} \ast \cdots \ast \Delta_{2r-1-3x_i,r} \ast \cdots \ast \Delta_{2r-1-3x_{d+1},r}\right) \le i\left( S \left(W_{r}^{\oplus (d+1)}\right) \right). \nonumber\end{eqnarray}

\vspace{0.2cm}

From Proposition $2.4(iv)$ we have $i\left( S \left(W_{r}^{\oplus (d+1)}\right) \right) = i\left( S^{(r-1) (d+1)-1}\right) = (r-1) (d+1)$.

\vspace{0.2cm}

We can calculate $\mbox{conn} (\Delta_{2r-1-3x_i,r})$, for all $i \in [d+1]$, using Theorem $2.13$, as follows.

\begin{eqnarray}\mbox{conn} (\Delta_{2r-1-3x_i,r})= \mbox{min} \left \{  2r-1-3x_i, r, \left \lfloor \frac{(2r-1-3x_i)+r +1}{3} \right \rfloor \right \}-2 = \nonumber \\ \mbox{min} \left \{  2r-1-3x_i, r, r-x_i \right \} \stackrel{(2x_i +1 \le r)}{=} (r-x_i)-2= r-x_i-2. \hspace{1cm} \nonumber \end{eqnarray}

\vspace{0.2cm}

It follows from Theorem $2.10$ that: \begin{eqnarray} \mbox{conn} (\Delta_{2r-1-3x_1,r} \ast \cdots \ast \Delta_{2r-1-3x_i,r} \ast \cdots \ast \Delta_{2r-1-3x_{d+1},r} ) \ge  \left [\displaystyle \sum_{i=1}^{d+1} (r-x_i-2) \right]  +2d = \nonumber \\ (d+1) (r-2) - \left( \displaystyle \sum_{i=1}^{d+1} x_i \right) +2d \stackrel{(\sum_{i=1}^{d+1} x_i \le d)}{\ge} (d+1)(r-1)-1. \hspace{1.5cm} \nonumber \end{eqnarray}

\vspace{0.2cm}

Consequently, from Theorem $2.12$ we have: \begin{eqnarray} i ( \Delta_{2r-1-3x_1,r} \ast \cdots \ast \Delta_{2r-1-3x_i,r} \ast \cdots \ast \Delta_{2r-1-3x_{d+1},r} ) \ge [(d+1)(r-1)-1]+2= \nonumber \\  (d+1)(r-1)+1 > (r-1) (d+1). \hspace{3cm}\nonumber \end{eqnarray}

\vspace{0.2cm}

Then, \begin{eqnarray}i(\Delta_{2r-1-3x_1,r} \ast \cdots \ast \Delta_{2r-1-3x_i,r} \ast \cdots \ast \Delta_{2r-1-3x_{d+1},r}) > (r-1) (d+1) =i \left(S\left(W_{r}^{\oplus(d+1)}\right)\right).\nonumber\end{eqnarray}

But it contradicts: \begin{eqnarray}i(\Delta_{2r-1-3x_1,r} \ast \cdots \ast \Delta_{2r-1-3x_i,r} \ast \cdots \ast \Delta_{2r-1-3x_{d+1},r}) \le i \left(S\left(W_{r}^{\oplus(d+1)}\right)\right). \nonumber \end{eqnarray} 

Therefore, there is no $(\mathbb{Z}_p)^n$-equivariant map \begin{eqnarray}\Delta_{2r-1-3x_1,r} \ast \cdots \ast \Delta_{2r-1-3x_i,r} \ast \cdots \ast \Delta_{2r-1-3x_{d+1},r} \longrightarrow S \left(W_{r}^{\oplus (d+1)}\right) \nonumber\end{eqnarray}.

\hspace{14.3cm} $\square$

\vspace{0.5cm}

\begin{example} \textbf{(I)} Note that if $x_1=0$, $x_2 = \cdots = x_{d+1}=1$, we have $\displaystyle \sum_{i=1}^{d+1} x_i =d \le d$ and $2x_i +1 \le r$, for all $i \in [d+1]$. This particular case of Theorem $5.1$ is exactly the stated in Theorem $4.2$.

\vspace{0.2cm}

\textbf{(II)} Let us consider $r=9$, $d=3$, $x_1=x_2=0$, $x_3=2$ and $x_4=1$. We have $\displaystyle \sum_{i=1}^{d+1} x_i =3 \le d=3$ and $2x_i +1 \le r$, for all $i \in [d+1]$. This is a case of the Theorem $5.1$, which does not occur in Theorem $4.2$, since $|C_3| = 2r-1-3x_3 = 11 < 14 =2r-4$.

\end{example}


\end{document}